\documentclass[10pt]{article}
\usepackage{amsmath,amsthm}
\usepackage{latexsym}
\usepackage{amssymb,amscd}
\usepackage{xspace}
\usepackage{amscd}
\usepackage{ graphicx, latexsym, amssymb, amscd, psfrag}
\usepackage{graphics}
\usepackage[frame,ps,dvips,matrix,arrow,curve,rotate]{xy}
\newtheorem{theo}{Theorem}[section]
\newtheorem{conj}[theo]{Conjecture}
\newtheorem{prop}[theo]{Proposition}

\newtheorem{lemm}[theo]{Lemma}
\newtheorem{coro}[theo]{Corollary}

\newcommand{\te}{\mathcal{T}} 
\newcommand{\vo}{\operatorname{Vol}}

\newcommand{\eLL}{L}

\newcommand{\Mod}{\operatorname{Mod}}

\newcommand{\vol}{\operatorname{Vol}}

\newcommand{\CC}{\beta}
\newcommand{\EE}{\mathbb E}
\newcommand{\ch}{\S}
\newcommand{\cC}{\mathcal{F}}
\newcommand{\PP}{\operatorname{Prob}}

\numberwithin{equation}{section}

\title{Growth of Weil-Petersson volumes and random hyperbolic surfaces of large genus}
\author{Maryam Mirzakhani\thanks{partially supported by an NSF grant.}} 

\begin{document}
\maketitle
\begin{section}{Introduction}
In this paper, we investigate the geometric properties of hyperbolic surfaces by studying the lengths of simple closed geodesics. 
The moduli space $\mathcal{M}_{g,n}$ of complete hyperbolic surfaces of genus $g \ge 2$ with $n$ punctures, 
is equipped with a natural notion of measure, which is induced by the {\em Weil-Petersson} symplectic form $\omega_{g,n}$ ($\ch \ref{BN}$). By a theorem of Wolpert, this form is the symplectic form of a K\"{a}hler noncomplete metric on the moduli space $\mathcal{M}_{g,n}$.
We describe the relationship between the behavior of lengths of simple closed geodesics on a hyperbolic surface 
and properties of the moduli space of such surfaces. 
First, we study the asymptotic behavior of Weil-Petersson volume $V_{g,n}$ of the moduli 
spaces of hyperbolic surfaces of genus $g$ with $n$ punctures as $g \rightarrow \infty.$
Then discuss some geometric properties of a random hyperbolic surface with respect to the Weil-Petersson measure as $g \rightarrow \infty$. \\
\noindent
{\bf Notation.} For any function $F : \mathcal{M}_{g} \rightarrow {\mathbb R}$, let 
$$\EE^g_{X\sim wp} (F(X))= \frac{\int_{\mathcal{M}_{g}} F(X) dX}{V_{g}},$$
where the integral is taken with respect to the Weil-Petersson volume form. 
 Also,
$$ \PP^g_{wp} (F(X) \leq C)= \EE^g_{X\sim wp} (G(X)), $$
 where $G(X)=1$ iff $F(X)\leq C$ and $G(X)=0$ otherwise. \\
 In this paper, $f_1(g) \asymp f_{2}(g)$ means that there exists a constant $C>0$ independent of $g$ such that $$\frac{1}{C} f_{2}(g) \leq f_1(g) \leq C f_{2}(g).$$
 Similarly, 
 $f_{1}(g)= O(f_{2}(g))$ means there exists a constant $C>0 $ independent of $g$ such that $$f_1(g) \leq C f_{2}(g).$$

\begin{subsection}{}
\noindent
{\bf Moduli spaces of hyperbolic surfaces with geodesic boundary components.}
The universal cover of $\mathcal{M}_{g,n}$ is the Teichm\" uller space $\te_{g,n}$. 
Every isotopy class of a closed curve on a hyperbolic surface $X \in \te_{g,n}$ contains a unique closed geodesic. 
Given a homotopy class of a closed curve $\alpha$ on a topological surface $S_{g,n}$ of genus $g$ with $n$ marked points and 
 $X \in \te_{g,n}$, let $\ell_{\alpha}(X)$ be the length of the unique 
geodesic in the homotopy class of $\alpha$ on $X$. This defines a length function $\ell_\alpha$ on the Teichm\"uller space $\te_{g,n}$.

When studying the behavior of hyperbolic length functions, it proves fruitful to consider more generally bordered hyperbolic surfaces with geodesic boundary components. Given $L=(L_{1},\ldots, L_{n})\in {\mathbb R}_{+}^{n}$, we consider the Teichm\"uller space $\te_{g,n}(L)$ of hyperbolic 
structures with geodesic boundary components of length $L_{1},\ldots, L_{n}$. Note that a geodesic of length zero is the same as a puncture. 
In fact, the space $\te_{g,n}(L)$ is naturally equipped with a symplectic form $\omega_{wp}$.
 The Weil-Petersson volume $V_{g,n}(L)$ of $\mathcal{M}_{g,n}(L_{1},\ldots,L_{n})$ is a polynomial 
in $L_{1}^2,\ldots, L_{n}^2$ of degree $3g-3+n$ and $V_{g,n}=V_{g,n}(0,\ldots,0)$.
We will show that in order to get bounds on the integrals of geometric functions over $\mathcal{M}_{g}$ we need to understand the asymptotics of the polynomial 
$V_{g,n}(L)$ as $g \rightarrow \infty.$
\end{subsection}
\begin{subsection}{}
\noindent
{\bf New results.}
Here we discuss the main results obtained in this paper:\\

\noindent
{\bf I):Asymptotic behavior of Weil-Petersson volumes.}
Peter Zograf has developed a fast algorithm for calculating the volume polynomials, and made several conjectures on the basis of the numerical data obtained by his algorithm \cite{Z:con}. 
\noindent
\begin{conj}[Zograf] For any fixed $n\geq 0$ 
$$V_{g,n} =(4\pi^2)^{2g+n-3} (2g-3+n)! \frac{1}{\sqrt{g \pi}} \left(1 + \frac{c_{n}}{g} + O\left(\frac{1}{g^{2}}\right)\right)$$
as $ g\rightarrow \infty.$ 
\end{conj}
Here $$V_{g,n}= \int\limits_{\overline{\mathcal{M}}_{g,n}} \omega_{g,n}^{3g-3+n}/ (3g-3+n)!.$$
\noindent 
In $\S \ref{Asym},$ we show : 
\begin{theo}\label{11}
For any $n \geq 0$: 
 $$\frac{V_{g,n+1}} {2gV_{g,n}}= 4\pi^{2}+O(\frac{1}{g}),$$
and
$$ \frac{V_{g,n}} {V_{g-1,n+2}}=1+O(\frac{1}{g})$$
as $g \rightarrow \infty.$
\end{theo}
These estimates imply that there exists $M>0$ such that 
\begin{equation}\label{compare}
 g^{-M} \leq \frac{V_{g,n}}{(4\pi^2)^{2g+n-3} (2g-3+n)! } \leq g^{M}. 
 \end{equation}

In order to prove this theorem, we discuss the asymptotics of all the coefficients of the volume polynomials $V_{g,n}(L)$ (see Theorem $\ref{relation1}$);
there are $$ \int\limits_{\overline{\mathcal{M}}_{g,n}} \!\!\psi_{1}^{d_{1}}\cdots \psi_{n}^{d_{n}}\cdot \omega^{3g-3+n-|{\bf d}|},$$
where for $1 \leq i \leq n,$ $\psi_{i} \in H^{2} (\mathcal{M}_{g,n}, {\mathbb Q})$ is the first Chern class of the tautological line bundle corresponding to the $i$-th puncture on 
$X \in \mathcal{M}_{g,n}$ (\ch \ref{BN}), and $|{\bf d}|=d_{1}+\ldots+d_{n}.$

 In \ch \ref{Asym} we apply known recursive formulas for these numbers and obtain some basic estimates for the intersection pairings 
of $\psi_{i}$ classes on $\overline{\mathcal{M}}_{g,n}$ as $g \rightarrow \infty$.\\

\noindent
{\bf II):Geometric behavior of surfaces of high genus.} In section \ch $\ref{Rs}$, we prove that as $g \rightarrow \infty$ the followings hold:
\begin{itemize}
\item
In \ch $\ref{sh}$, we show that the probability that a random Riemann surface has a short non-separating simple
closed geodesic is asymptotically positive. More precisely, 
let $\ell_{sys}(X)$ denote the length of the shortest simple closed geodesic on $X$.
Then for any small (but fixed) $\epsilon>0$, as $g \rightarrow \infty$
$$\PP^g_{wp} (\ell_{sys}(X)< \epsilon) \asymp \epsilon^2.$$

 \item 
However, {\it separating} simple closed geodesics tend to be much longer \ch $\ref{ss}$.
Let $\ell^{s}_{sys}(X)$ denote the length of the shortest {\it separating} simple closed geodesic on $X$. 
we show that 
$$ \PP^g_{wp} (\ell^{s}_{sys}(X)< m_{1} \log(g))=O(\log(g) g^{(m_{1}/2-1)})),$$
and
$$\EE^g_{X\sim wp} (\ell^{s}_{sys}(X)) \asymp \log(g)$$
as $g \rightarrow \infty.$
\item 
 Similarly, using the asymptotics of $V_{g,n}(L)$ we get bounds for the expected length of the shortest simple closed geodesic of a given combinatorial type. In particular, the shortest simple closed geodesic separating the surface into two roughly equal areas has length at least linear in $g$.
Moreover, in \ch $\ref{ch}$ we show that the Cheeger constant $h(X)$ of a random Riemann surfaces $X \in \mathcal{M}_{g}$ is bounded from below by a universal constant .
More precisely, as $g \rightarrow \infty$ 
$$\PP^g_{wp}\left (h(X) \leq \frac{\ln(2)}{\pi+\ln(2)}\right) \rightarrow 0.$$
By Cheeger's theoem the smallest positive eigenvalue of the Laplacian on a generic point $X$ is $\geq \frac{1}{4} C_{h}^{2}$,
where $C_{h}= \frac{\ln(2)}{\pi+\ln(2)}.$ 
\item Finally, we show that a generic hyperbolic surface in $\mathcal{M}_{g}$ has a small diameter, with a large embedded ball \ch $\ref{di}$. 
More precisely, as $g \rightarrow \infty$
$$\PP^g_{wp}(\operatorname{diam}(X) \geq C_{d}\log(g)) \rightarrow 0,$$
and
$$\EE^g_{X\sim wp} (\sqrt{\operatorname{diam}(X)}) \asymp \sqrt{\log(g)}.$$
Also,
$$\PP^g_{wp}(\operatorname{Emb}(X) \leq C_{E}\log(g)) \rightarrow 0,$$
and
$$\EE^g_{X\sim wp} (\operatorname{Emb}(X)) \asymp \log(g)$$
where $\operatorname{Emb}(X)$ is the radius of the largest embedded ball in $X.$
Here $C_{E}=\frac{1}{3}, $ and $C_{d}=5.$
\end{itemize}
We remark that none of the costants in these statements are sharp. However, in this paper for the sake of simplicity, we chose the simpler proofs which would give weaker constants.\\
%{\bf Notation.} In this note,\noindent
\end{subsection}
\begin{subsection}{}\label{O}
Our main tool is the close relationship between the Weil-Petersson geometry of $\mathcal{M}_{g,n}$ and the lengths of simple closed geodesics on surfaces in $\mathcal{M}_{g}$. Here we discuss one application of this relationship.

Let $\mathcal{S}_{g,n}$ denote the set of homotopy classes of non-trivial simple closed curves on a topological surface $S_{g,n}$ of genus $g$ with $n$ marked points.
For any $\gamma \in \mathcal{S}_{g,n},$ let $S_{g,n}-\gamma$ denote the surface obtained by cutting the surface $S_{g,n}$ 
along $\gamma.$
Given $\alpha_{1},\alpha_{2} \in \mathcal{S}_{g}$, we say $\alpha_{1} \sim \alpha_{2}$ if $\alpha_{1}$ and $\alpha_{2}$ are of the same {\it type}; 
that is $S_{g}-\alpha_{1}$ is homeomorphic to $S_{g}-\alpha_{2}.$
Given a connected simple closed curve $\gamma \in \mathcal{S}_{g}$ and $f: {\mathbb R}_{+} \rightarrow {\mathbb R}_{+}$ 
define 
 $f_{\gamma} : \te_{g} \rightarrow {\mathbb R}_{+}$ by
$$f_{\gamma} (X)= \sum_{\alpha \sim \gamma} f(\ell_{\alpha}(X)).$$
We have
\begin{equation}\label{lemineq}
\int_{{\mathcal M}_{g}}\!\! f_{\gamma}(X)\;dX=\int\limits_{0}^{\infty}
f(t)\;t\;V_{g,2}(t,t)\; dt.
\end{equation}
For the general case of this formula see Theorem $\ref{integrate}.$
The main idea is that the decomposition of the surface along $\gamma$ gives rise
to a description of a cover of $\mathcal{M}_{g}$ in terms
of moduli spaces corresponding to simpler surfaces. \\
\end{subsection}
\begin{subsection}{}
\noindent
{\bf Notes and remarks.}
\begin{itemize}
\item A recursive formula for the Weil-Petersson volume of the moduli space 
of punctured spheres was obtained by Zograf \cite{Z:p}. Moreover, Zograf and 
Manin have obtained generating functions for the Weil-Petersson volume of 
$\mathcal{M}_{g,n}$\cite{MZ}. See also \cite{KMZ:w}.
The following exact asymptotic formula was proved
in \cite{MZ}.
\begin{theo}\label{z}
 There exists $C>0$ such that for any fixed $g\geq 0$ 
\begin{equation}\label{zn} 
V_{g,n} = n! C^{n} n^{(5g-7)/2} (a_{g} + O (1/n)),
\end{equation}
as $ n \rightarrow \infty.$
\end{theo}
 \item In \cite{Gr:vol}, it is shown that for a fixed $n>0$ there are $c_{1},c_{2}>0$ such that 
 $$c_{2}^g (2g)! < \vo(\mathcal{M}_{g,n}) < c_{1}^g (2g)!.$$ 
 This result was extended to the case of $n=0$ in \cite{ST:vol}. 
 Note that these estimates do not give any information about the growth of $$B_{g,n}=V_{g,n}/V_{g-1, n+2}$$ and $$C_{g,n}=V_{g,n+1}/(2g V_{g,n})$$ when 
 $g \rightarrow \infty.$
 \item Penner has developed a different method for calculating the Weil-Petersson volume of the moduli spaces of curves with marked points by using decorated Teichm\"uller theory \cite{P:vol}.
 \item In \cite{B:M} Brooks, and Makover developed a method for the study of {\it typical} Riemann surfaces with large genus by 
using trivalent graphs. In this model the expected value of the systole of a random Riemann surface turns 
out to be bounded (independent of the genus) \cite{M:M}. See also \cite{Ga}.
We will see in this note that a random Riemann surface with respect to the Weil-Petersson volume form has similar features. 
However, it is not clear how the measure induced by their model is related to the measure induced by the Weil-Petersson volume.
\item The distribution of hyperbolic surfaces of genus $g$ produced
randomly by gluing Riemann surfaces with long geodesic boundary components is closely related to the volume form induced by $\omega$ on 
$\mathcal{M}_{g,n}$. See \cite{M:AB} for details.
 \end{itemize}
 \end{subsection}
 \noindent
 {\bf Questions.}
\begin{itemize}
\item 
It would be useful to know the asymptotics of $$\frac{V_{g,n(g)}}{V_{g-1,n(g)+2}},$$
where $n(g) \rightarrow \infty$ as $g \rightarrow \infty$.
 Note that by Theorem $\ref{z}$ and Theorem $\ref{11}$, we know the asymptotics of $V_{g}/V_{g-1,2}$ and $V_{1,2g-4}/V_{0,2g-2}$. However, we don't know 
 much about the behavior of the sequence
 $$ V_{g}\;,\; V_{g-1,2}\;,\; \ldots\;, \;V_{0,2g}$$
 as $g \rightarrow \infty.$
 
 \item As in Theorem \ref{relation1}, when $n=1$ the volume polynomial can be written as
 $$V_{g,1}(L)= \sum_{k=0}^{3g-2} \frac{a_{g,k}}{(2k+1)!} L^{2k}, $$
 where $a_{g,k}$ are rational multiples of powers of $\pi.$
 It would be helpful to understand the asymptotics of $a_{g,k}/a_{g,k+1}$ for an arbitrary $k$ (which can grow with $g$).
 Note that $a_{g,0}=V_{g,1}.$ In Theorem $\ref{first}$(a), we show that for given $i\geq 0$
 $$\lim_{g \rightarrow \infty} \frac{a_{g,i+1}}{a_{g,i}}=1.$$
 On the other hand, it is known that \cite{I:Z}
 $$ \int\limits_{\overline{\mathcal{M}}_{g,}} \!\!\psi_{1}^{3g-2}=\frac{1}{24^{g}g!},$$
 and hence
 $$\frac{a_{g,3g-2}}{a_{g,0}}\rightarrow 0,$$
 as $g \rightarrow \infty.$
 \item The results obtained in this paper are only small steps towards understanding the geometry of random hyperbolic surfaces 
 of large genus. Many interesting questions about such random surfaces are open.
 Investigating geometric properties of random Riemann surfaces could shed some light on the asymptotics geometry of $\mathcal{M}_{g}$ as $g \rightarrow \infty,$ see \cite{C:P}, \cite{T}, and \cite{Hu} for some results in this direction. 
 \end{itemize}
 
 \noindent
 {\bf Acknowledgement.} 
 I would like to thank P. Zograf for many helpful and illuminating discussions regarding the growth of Weil-Petersson volumes. 
 I am also grateful to Rick Schoen and Jan Vondrak.
 %%WHAT???????????????????????????????????????????????????????????????????
\end{section}
\begin{section}{Background and notation}\label{BN}
In this section, we recall definitions and known results about the geometry 
of hyperbolic surfaces and properties of their moduli spaces. For more details see \cite{M:In}, \cite{B:book} and \cite{W:M}.

\begin{subsection}{}
\noindent
{\bf Teichm\"uller Space.} A point in the {\it Teichm\"uller space} $\te(S)$
is a complete hyperbolic surface $X$ equipped with a diffeomorphism $f: S
\rightarrow X$. The map $f$ provides a {\it marking} on $X$ by $S$. Two
marked surfaces $f:\; S \rightarrow X$ and $g:\; S \rightarrow Y$
define the same point in $\te(S)$ if and only if $f\circ g^{-1} : Y \rightarrow X$
is isotopic to a conformal map. When $\partial S$ is nonempty, consider 
hyperbolic Riemann surfaces homeomorphic to $S$ with geodesic
boundary components of fixed length. Let $A=\partial S $ and $\eLL=(\eLL_{\alpha})_{\alpha \in A} \in {\mathbb R}_{+}^{|A|}$. A point $X \in\te_{g,n}(\eLL)$ is a marked hyperbolic surface 
with geodesic boundary components such
that for each boundary component $\beta \in \partial S$, we have 
 $$\ell_{\beta}(X)=\eLL_{\beta}.$$ 
 By convention, a geodesic of length zero is a cusp and we have 
$$\te_{g,n}=\te_{g,n}(0,\ldots,0).$$
Let $\operatorname{Mod}(S)$
denote the mapping class group of $S$, or the
group of isotopy classes of orientation preserving self
homeomorphisms of $S$ leaving each boundary component setwise fixed. The mapping
class 
group 
$\operatorname{Mod}_{g,n}=\operatorname{Mod}(S_{g,n})$ acts 
on $\te_{g,n}(\eLL)$ by changing the marking. The quotient space
$$\mathcal{M}_{g,n}(\eLL)=\mathcal{M}(S_{g,n},\ell_{\CC_{i}}=\eLL_{i})=
\te_{g,n}(\eLL_{1},\ldots,\eLL_{n})/\operatorname{Mod}_{g,n}$$
is the moduli space of Riemann surfaces homeomorphic to $S_{g,n}$ 
with $n$ boundary components of length $\ell_{\CC_{i}}=\eLL_{i}$.
Also, we have 
$$\mathcal{M}_{g,n}=\mathcal{M}_{g,n}(0,\ldots,0).$$
By work of Goldman \cite{G:S},
the space $\te_{g,n}(\eLL_{1},\ldots,\eLL_{n})$ carries a natural symplectic
form invariant under the action of the mapping class group.
This symplectic form is called the
{\it Weil-Petersson symplectic form,} and denoted by $\omega$ or $\omega_{wp}$.
This symplectic form is in fact the K\"ahler form of a K\"ahler metric \cite{IT:book}.\\
By work of Wolpert, over Teichm\"uller space the
Weil-Petersson symplectic structure has a simple form in Fenchel-Nielsen
coordinates \cite{W:F}. \\
\noindent
{\bf The Fenchel-Nielsen coordinates.} A {\it pants decomposition} 
of $S$ is a set of disjoint simple closed curves which decompose 
the surface into pairs of pants. Fix a system of 
pants decomposition of $S_{g,n}$, 
$\mathcal{P}=\{\alpha_{i}\}_{i=1}^{k}$, where
$k=3g-3+n$. For a marked hyperbolic surface 
$X \in \te_{g,n}(\eLL)$, the {\it Fenchel-Nielsen coordinates} 
associated with 
$\mathcal{P}$, $\{\ell_{\alpha_{1}}(X),\ldots,\ell_{\alpha_{k}}(X),
\tau_{\alpha_{1}}(X),\ldots,
\tau_{\alpha_{k}}(X)\}$, 
 consists of the 
set of lengths of 
all geodesics used in the decomposition and the set of the 
{\it twisting} parameters used to glue the pieces.
We have an isomorphism 
$$\te_{g,n}(\eLL) \cong {\mathbb R}_{+}^{\mathcal{P}}\times {\mathbb R}^{\mathcal{P}} $$ 
by the map
$$X \rightarrow (\ell_{\alpha_{i}}(X),\tau_{\alpha_{i}}(X)).$$
See \cite{B:book} for more details.
\begin{theo}[Wolpert]\label{wolp}
The Weil-Petersson symplectic form is given by
$$\omega_{wp} =\sum\limits_{i=1}^{k} d\ell_{\alpha_{i}}\wedge
d\tau_{\alpha_{i}}.$$
\end{theo}
By Theorem \ref{wolp} the natural twisting around $\alpha$ 
is the Hamiltonian flow of the length function of $\alpha$.\\ 
\noindent

\end{subsection}
\begin{subsection}{}
\noindent
{\bf Integrating geometric functions over moduli spaces.}
Here, we discuss a method for integrating certain geometric functions over $\mathcal{M}_{g,n}(\eLL)$.
 Let $Y \in \te_{g,n}.$ For a simple closed curve $\gamma$ on $S_{g,n}$, let
$[\gamma]$ denote the homotopy class of $\gamma$ and let 
$\ell_{\gamma}(Y)$ denote the hyperbolic length of the geodesic 
representative of $[\gamma]$ on $Y$. 
To each simple closed curve
$\gamma$ on $S_{g,n}$, we associate the set 
$$\mathcal{O}_{\gamma}=\{[\alpha]\;| \alpha \in
\operatorname{Mod}_{g,n} \cdot \gamma\}$$ 
of homotopy classes of simple closed curves in the $\operatorname{Mod}_{g,n}$-orbit of $\gamma$
on $X \in \mathcal{M}_{g,n}$. 
Given a function $f:{\mathbb R}_{+} \rightarrow {\mathbb R}_{+}$, and a multicurve $\gamma$ on $S_{g,n}$
define $$f_{\gamma}: \mathcal{M}_{g,n} \rightarrow {\mathbb R}$$
by 
\begin{equation}\label{function}
f_{\gamma}(X)=\sum\limits_{[\alpha] \in \mathcal{O}_{\gamma}} 
f(\ell_{\alpha}(X)).
\end{equation}
Let $S_{g,n}(\gamma)$ be the result of cutting the surface $S_{g,n}$ 
along $\gamma$; that is $S_{g,n}(\gamma) \cong S_{g,n}-U_{\gamma}$, where $U_{\gamma}$ is 
an open neighborhood of $\gamma$ homeomorphic to $\gamma \times (0,1)$.
Thus $S_{g,n}(\gamma)$ is a possibly disconnected compact surface 
with $n+2$ boundary components. We define
$\mathcal{M}(S_{g,n}(\gamma),\ell_{\gamma}=t)$ to be the moduli space 
of Riemann surfaces homeomorphic to $S_{g,n}(\gamma)$ such that the
lengths of the 2 boundary components corresponding to $\gamma$ are equal 
to $t$. Then we have (\cite{M:In}):
\begin{theo}\label{integrate}
For any multicurve $\gamma=\sum\limits_{i=1}^{k} c_{i} \gamma_{i}$, the integral of $f_{\gamma}$ over $\mathcal{M}_{g,n}(\eLL)$ with 
respect to the Weil-Petersson volume form is given 
by 
$$\int \limits_{\mathcal{M}_{g,n}(\eLL)}\!\! 
f_{\gamma}(X)\,
dX=\frac{2^{-M(\gamma)}}{|\operatorname{Sym}(\gamma)|}\;\int\limits_{{\bf x}
 \in {\mathbb R}_{+}^{k}} f(|{\bf x}|)\; V_{g,n}(\Gamma,{\bf x},\CC,\eLL)\; {\bf
x}\cdot d {\bf x},
$$
where $\Gamma=(\gamma_{1},\ldots,\gamma_{k})$, $|{\bf x}|= \sum\limits_{i=1}^{k} c_{i}\; x_{i}$, 
$ {\bf x}\cdot d{\bf x}=x_{1} \cdots
x_{k} \cdot dx_{1}\wedge \cdots \wedge dx_{k}$,
and 
$$M(\gamma)=| \{i | \gamma_{i} \mbox{ \; \mbox{separates} off a one-handle from \;} S_{g,n} \}|.$$
\end{theo}
Given a multicurve $\gamma=\sum_{i=1}^{k} c_{i}\gamma_{i}$, the symmetry group of $\gamma$, 
$\operatorname{Sym}(\gamma)$, is defined by 

$$\operatorname{Sym}(\gamma)= \operatorname{Stab}(\gamma)/\cap_{i=1}^{k} \operatorname{Stab}(\gamma_{i}). $$

Recall that given ${\bf x}=(x_{1},\ldots,x_{k}) \in {\mathbb R}_{+}^{k}$, $V_{g,n}(\Gamma,{\bf x},\CC,\eLL)$ is defined by 
 $$V_{g,n}(\Gamma,{\bf x},\CC,\eLL)= \vo(\mathcal{M}(S_{g,n}(\gamma),
\ell_{\Gamma}={\bf x},\ell_{\CC}=\eLL)).$$ 
Also,
$$V_{g,n}(\Gamma,{\bf x},\CC,\eLL)=\prod \limits_{i=1}^{s}
V_{g_{i},n_{i}}(\ell_{A_{i}}),$$
where 
\begin{equation}\label{disj}
S_{g,n}(\gamma)=\bigcup_{i=1}^{s} S_{i}\;,
\end{equation}
$S_{i} \cong S_{g_{i},n_{i}},$ and 
$\; A_{i}= \partial
S_{i}.$ \\
 By Theorem $\ref{integrate}$ integrating $f_{\gamma}$, even for a compact Riemann surface, reduces
to the calculation of 
volumes of moduli spaces of bordered Riemann surfaces.
%This formula can be used to relate the growth of the number of simple closed geodesics on $X \in \mathcal{M}_{g}$ to the volume polynomials \cite{M:C}.\\
%See $\ch \ref{ground}$ for the definition of $\operatorname{Sym}(\gamma)$, the
 %symmetry group of $\gamma=\sum_{i=1}^{k} c_{i} \gamma_{i}$.\\
\noindent
%%%%%%%%%%%%MORE%%%%%%%%%%%%
{\bf Remark.} %The terms $|\operatorname{Sym}(\gamma)|$ and $M(\gamma)$
%are not equal to one when the multicurve $\gamma$ has some extra symmetry.
%Note that $|\operatorname{Sym}(\gamma)| \not = 1$ will impose a nontrivial restriction on the $c_{i}$'s.
 Let $g \in \operatorname{Sym}(\gamma)$, where $\gamma=\sum_{i=1}^{k} c_{i} \gamma_{i}$. Then 
 $g(\gamma_{i})=\gamma_{j}$ implies that $c_{i}=c_{j}$. 
\end{subsection}
\begin{subsection}{}
\noindent
{\bf Connection with the intersection pairings of tautological line bundles.}
The moduli space $\mathcal{M}_{g,n}$ is endowed with natural cohomology classes. 
 When $n>0,$ there are $n$ tautological line bundles defined on 
$\overline{\mathcal{M}}_{g,n}$ as follows.
We can define $\mathcal{L}_{i}$ in the orbifold sense whose fiber at the point 
$(C,x_{1},\ldots,x_{n})\in \overline{\mathcal{M}}_{g,n}$
is the cotangent space of $C$ at $x_{i}.$ 
Then $\psi_{i}= c_{1}(\mathcal{L}_{i}) \in H_{2} (\overline{\mathcal{M}}_{g,n}, {\mathbb Q})$. 
Note that although the complex curve $C$ may have nodes, $x_{i}$ never
coincides with the singular points.
 See \cite{Harris:book} and \cite{AC} for more details.

In \cite{M:JAMS}, we use the symplectic geometry of moduli spaces of bordered Riemann surfaces to relate these intersection pairings to 
the volume polynomials. This method allows us to read off the intersection numbers of tautological line bundles from the volume polynomials:
\begin{theo}\label{relation1}
In terms of the above notation,
$$\vo(\mathcal{M}_{g,n}(\eLL_{1},\ldots,\eLL_{n}))=\sum\limits_{ |{\bf d}|
 \leq 3g-3+n}C_{g}({\bf d})\;\; \eLL_{1}^{2d_{1}}\ldots \eLL_{n}^{2d_{n}},$$ 
where ${\bf d}=(d_{1},\ldots, d_{n}),$ and $C_{g}({\bf d})$ is equal to
$$ \frac{2^{m(g,n) |{\bf d}|}} {2^{|{\bf d}|} \; |{\bf d}|!\;(3g-3+n-|{\bf d}|)! }\; \int\limits_{\overline{\mathcal{M}}_{g,n}} \!\!\psi_{1}^{d_{1}}\cdots \psi_{n}^{d_{n}}\cdot \omega^{3g-3+n-|{\bf d}|}.$$ 
Here $ m(g,n)= \delta(g-1) \times \delta(n-1)$,
${\bf d}!= \prod_{i=1}^{n} d_{i}!$, and $|{\bf d}|=\sum_{i=1}^{n} d_{i}.$
\end{theo}
\noindent
{\bf Remark.}
 We warn the reader that there are some small differences in the normalization of the Weil-Petersson volume 
form in the literature; in this paper, 
$$V_{g,n}=V_{g,n}(0,\ldots,0)=\frac{1}{(3g-3+n)!} \int_{\mathcal{M}_{g,n}} \omega^{3g-3+n}$$
which is slightly different from the notation used in \cite{Z:con} and \cite{ST:vol}.
Also, in \cite{Z:p} the Weil-Petersson K\"ahler form is $1/2$ the 
imaginary part of the Weil-Petersson pairing, while here the factor $1/2$ does not
appear. So our answers are different by a power of $2$. 

\end{subsection}
\end{section}

\begin{section}{Asymptotic behavior of Weil-Petersson volumes}\label{Asym} 
%This formula can be used to relate the growth of the number of simple closed geodesics on $X \in \mathcal{M}_{g}$ to the volume polynomials \cite{M:C}.\\
%See $\ch \ref{ground}$ for the definition of $\operatorname{Sym}(\gamma)$, the
 %symmetry group of $\gamma=\sum_{i=1}^{k} c_{i} \gamma_{i}$.\\
In this section, we study the asymptotics of $V_{g,n}({\bf L})=\vo(\mathcal{M}_{g,n}(L_{1}\ldots,L_{n}))$ as $g \rightarrow \infty.$ \\
\noindent
{\bf Notation.} For ${\bf d}=(d_{1},\ldots,d_{n})$ with $ d_{i} \in {\mathbb N} \cup \{0\}$ and $|{\bf d}|=d_{1}+\ldots+d_{n} \leq 3g-3+n,$ let $d_{0}= 3g-3-|{\bf d}|$ and define
$$ [ \prod_{i=1}^{n} \tau_{d_{i}}]_{g,n}= \frac{\prod_{i=1}^{n} (2d_{i}+1)! 2^{|{\bf d}|}} {\prod_{i=0}^{n} d_{i}!} \int_{\overline{\mathcal{M}}_{g,n}} \psi_{1}^{d_{1}} \cdots \psi_{n}^{d_{n}} 
\omega^{d_{0}}=$$
$$ =\frac{\prod_{i=1}^{n} (2d_{i}+1)!! 2^{2|{\bf d}|} (2\pi^2)^{d_{0}} } {d_{0}!} \int_{\overline{\mathcal{M}}_{g,n}} \psi_{1}^{d_{1}} \cdots \psi_{n}^{d_{n}} \kappa_{1}^{d_{0}},$$
 where $\kappa_{1}= \omega/(2 \pi^2)$ is the first Mumford class on $\overline{\mathcal{M}}_{g,n}$ \cite{AC}. By Theorem $\ref{relation1}$ for $L=(L_{1},\ldots,L_{n})$ we have:
\begin{equation}\label{re2}
V_{g,n}(2L)=\sum_{|{\bf d}| \leq 3g-3+n} [\tau_{d_{1}},\ldots \tau_{d_{n}}]_{g,n} \; \frac{L_{1}^{2d_{1}}}{(2d_{1}+1)!}\cdots \frac{L_{n}^{2 d_{n}}}{(2d_{n}+1)!}.
\end{equation} 
\begin{subsection}{}
\noindent
{Recursive formulas for the intersection pairings.} 
 Given ${\bf d}=(d_{1},\ldots, d_{n})$ with $|{\bf d}| \leq 3g-3+n,$ the following recursive formulas hold:\\
\noindent
{\bf I.}
$$
 [ \tau_{0} \tau_{1} \prod_{i=1}^{n} \tau_{d_{i}}]_{g,n+2} = [\tau_{0}^{4} \prod_{i=1}^{n} \tau_{d_{i}}]_{g-1,n+4}+$$
$$+\frac{1}{2}\sum_{g_{1}+g_{2}=g\; \atop \{1,\ldots, n\}=I\amalg J} [ \tau_{0}^{2} \; \prod_{i\in I } \tau_{d_{i}}]_{g_{1},|I|+2} \cdot [\tau_{0}^{2}\prod_{i\in J} \tau_{d_{i}}]_{g_{2},|J|+2},$$

\noindent
{\bf II.}
$$
(2g-2+n) [ \prod_{i=1}^{n} \tau_{d_{i}}]_{g,n}= \frac{1}{2} \sum_{L=0}^{3g-3+n} (-1)^{L} (L+1) \frac{\pi^{2L}}{(2L+3)!} [ \tau_{L+1} \; \prod_{i=1}^{n} \tau_{d_{i}}]_{g,n+1}.
$$

\noindent
{\bf III.}
Let $a_{0}=1/2,$ and for $n\geq 1$, $$a_{n} =\zeta(2n) (1-2^{1-2n}).$$
Then we have
$$
[ \tau_{d_{1}},\ldots,\tau_{d_{n}}]_{g,n}=\sum_{j=2}^{n}\mathcal{A}^{j}_{{\bf d}} +\frac{1}{2} \mathcal{B}_{{\bf d}}+ \frac{1}{2} \mathcal{C}_{{\bf d}},
$$
where 
\begin{equation}\label{A}
 \mathcal{A}^{j}_{{\bf d}}= \sum_{L=0}^{d_{0}} (2d_{j}+1) \; a_{L} [\tau_{d_{1}+d_{j}+L-1}, \prod_{i\not=1,j} \tau_{d_{i}}]_{g,n-1},\end{equation}
 
\begin{equation}\label{B}
 \mathcal{B}_{{\bf d}}= \;\sum_{L=0}^{d_{0}} \sum_{k_{1}+k_{2}=L+d_{1}-2} a_{L} [\tau_{k_{1}} \tau_{k_{2}} \prod_{i\not=1} \tau_{d_{i}}]_{g-1,n+1},
 \end{equation}
 and 
\begin{equation}\label{C}
\mathcal{C}_{{\bf d}}= \sum_{I \amalg J=\{2,\ldots,n\}\atop 0\leq g' \leq g} \sum_{L=0}^{d_{0}} \sum_{k_{1}+k_{2}=L+d_{1}-2} a_{L} \; [\tau_{k_{1}} \prod_{i\in I } \tau_{d_{i}}]_{g',|I|+1} \times [ \tau_{k_{2}} \prod_{i\in J} \tau_{d_{i}}]_{g-g',|J|+1}.
\end{equation}

\noindent
{\bf References}. 
\begin{itemize}

\item For results on the relationship between the Weil-Petersson volumes and the intersections of $\psi-$classes on $\overline{\mathcal{M}}_{g,n}$ see \cite{W} and \cite{AC}. 
An explicit formula for the volumes in terms of the intersection of $\psi-$classes was developed in \cite{KMZ:w}.

\item Formula $({\bf I})$ is a special case of Proposition $3.3$ in \cite{LX:higher}.

\item For different proofs of $({\bf II})$ see \cite{DN:cone} and 
\cite{LX:higher}. The proof presented in \cite{DN:cone} uses 
the properties of moduli spaces of of hyperbolic surfaces with cone points. 

 \item For a proof of $({\bf III})$ see \cite{M:In}; in view of Theorem \ref{relation1}, $({\bf III})$ can be interpreted as a recursive formula for the volume of $\mathcal{M}_{g,n}(\eLL)$
in terms of volumes of moduli spaces of Riemann surfaces that we get 
by removing a pair of pants containing at least one boundary
component of $S_{g,n}$. See also \cite{M:S} and \cite{LX:M}.

 \item If $d_{1}+\ldots+d_{n}=3g-3+n,$ $({\bf III})$ gives rise to a recursive formula for the intersection pairings of $\psi_{i}$ classes
 which is the same as the Virasoro
constraints for a point. See also \cite{MuS}. 
For different proofs and discussions related to these relations see \cite{W}, \cite{Ko:int}, \cite{OP}, \cite{M:JAMS}, \cite{KL:W}, and \cite{EO:WK}.
\end{itemize}
\noindent
{\bf Remarks.}
\begin{itemize}
\item In terms of the volume polynomials equation
${(\bf II)}$ can be written as (\cite{DN:cone}):
$$\frac{\partial V_{g,n+1}}{\partial L}(L, 2\pi i) = 2\pi i (2g-2+ n)V_{g,n}(L).$$
When $n=0,$ 
 $$V_{g,1}(2\pi i)=0,$$
 and
\begin{equation}\label{zero}
\frac{\partial V_{g,1}}{\partial L} (2\pi i) = 2\pi i (2g-2)V_{g}.
\end{equation}

\item 
Note that $({\bf III})$ applies only when $n>0.$ In the case of $n=0,$ $(\ref{zero})$ allows us to prove necessary estimates for the growth of $V_{g,0}.$ 

\item Although $({\bf III})$ has been described in purely
combinatorial terms, it is closely related to the topology 
of different types of pairs of pants in a surface. 

 \item In this paper, we are mainly interested in the intersection parings only containing $\kappa_{1}$ and $\psi_{i}$ classes. For generalizations of $({\bf III})$ to the case of higher Mumford's $\kappa$ classes see \cite{LX:higher} and \cite{E:M}.

\item We will show that $n$ is fixed and $g \rightarrow \infty$ both terms $\mathcal{A}_{{\bf d}},$ and $\mathcal{B}_{{\bf d}}$ in $({\bf III})$ contribute to $V_{g,n}=[\tau_{0},\ldots, \tau_{0}]_{g}$. More precisely, for 
${\bf d}=(0,\ldots,0)$ 
 $$\frac{\mathcal{B}_{{\bf d}}}{\mathcal{A}_{{\bf d}}} \asymp 1.$$
On the other hand, for ${\bf d}=(0,\ldots, 0)$ the contribution of $\mathcal{C}_{\bf d}$ in {\bf III} is negligible. More precisely, we will see that
$\frac{\mathcal{C}_{{\bf d}}}{V_{g,n}}=O(1/g).$
\end{itemize}

\end{subsection}
\begin{subsection}{}
\noindent
{\bf Basic estimates for the intersection pairings.}
The main advantage of using $({\bf III})$ is that all the coefficients are positive. Moreover, it is easy to check that
$$ a_{n}=\zeta(2n) (1-2^{1-2n})=\frac{1}{(2n-1)!} \int_{0}^{\infty} \frac{t^{2n-1}}{1+e^{t}}\; dt.$$
Hence,
$$a_{n+1}-a_{n}=\int_{0}^{\infty} \frac{1}{(1+e^{t})^{2}} \left(\frac{t^{2n+1}}{(2n+1)!}+\frac{t^{2n}}{2n!}\right) dt.$$
As a result, we have:
\begin{lemm}\label{aa}
In terms of the above notation, $\{a_{n}\}_{n=1}^{\infty}$ is an increasing sequence. 
Moreover, 
$\lim_{n\rightarrow \infty} a_{n}=1,$
and
\begin{equation}\label{abound}
a_{n+1}-a_{n} \asymp 1/2^{2n}.
\end{equation}
\end{lemm}
Using this observation and $(\ref{re2})$ one can prove the following general estimates:

 \begin{lemm}\label{obser}
 In terms of the above notation, the following estimates hold:
 \begin{enumerate}
 \item $$[\tau_{d_{1}},\tau_{0},\ldots, \tau_{0}]_{g,n} \leq [\tau_{0},\ldots, \tau_{0}]_{g,n}=V_{g,n},$$
and in case of ${\bf d}=(1,0,\ldots,0)$ we have
$$ [\tau_{1},\tau_{0}, \ldots, \tau_{0}]_{g,n} \asymp [\tau_{0},\ldots, \tau_{0}]_{g,n}.$$
\item More generally,
$$[\tau_{d_{1}},\ldots, \tau_{d_{n}}]_{g,n} \leq (2d_{1}+1)\cdots (2d_{n}+1) V_{g,n}, $$
and 
\begin{equation}\label{upper}
V_{g,n}(2L_{1},\ldots, 2L_{n}) \leq e^{L} V_{g,n},
\end{equation}
where $L=L_{1}+\ldots+L_{n}.$
\item 
for any $g,n \geq 0$, 
\begin{equation}\label{yekido}
V_{g,n+2} \geq V_{g-1,n+4},\;\; \mbox{and}\;\; \; V_{g,n+1} > \frac{(2g-2+n)}{b} V_{g,n},
\end{equation}\\
where $b=\sum_{L=0}^{\infty} \pi^{2L} (L+1) /2(2L+3)!.$ \\
\end{enumerate}
\end{lemm}
\noindent
{\bf Proof.}
Part $(1)$ and $(2)$ follow by comparing the contributions of $\mathcal{A}_{\bf d},$ $\mathcal{B}_{\bf d}$, and $\mathcal{C}_{\bf d}$
 for $(d_{1},d_{2},\ldots,d_{n})$, $(d_{1},0,\ldots,0)$ and $(0,\ldots,0)$ in $({\bf III})$. See $(\ref{A})$, $(\ref{B}),$ and $(\ref{C}.)$
 Then $(\ref{re2})$ implies ($\ref{upper}$).
 
Moreover, since $$[\tau_{1},\tau_{0},\ldots,\tau_{0}] \leq V_{g,n}, \;\;\mbox{and}\;\; [\tau_{0},\ldots,\tau_{0}]_{g'} \geq 0 $$
equation $({\bf I})$ for ${\bf d}=(0)$ implies that for any $n\geq 0$, $V_{g,n+2} \geq V_{g-1,n+4}$. Similarly, since $$[\tau_{L+1},\tau_{0},\ldots,\tau_{0}] \leq V_{g,n+1},$$
equation $({\bf II})$ for ${\bf d}=0$ implies $b V_{g,n+1} > 2(2g-2+n) V_{g,n}.$ 
\hfill $\Box$\\

\noindent
{\bf Remarks.} 
\begin{itemize}
\item 
A stronger lower bound for $\frac{V_{g,n+1}}{(2g-2+n) V_{g,n}} $ was obtained in \cite{ST:vol}. But in this paper, we will use only $(\ref{yekido}).$

\item We will show that as $g \rightarrow \infty$ the first inequality of $(\ref{yekido})$ is asymptotically 
sharp. However, $(\ref{zn})$ implies that when $g$ is fixed and $n$ is large this inequality is far from being sharp; in fact, given $g \geq 1$ as $n \rightarrow \infty$ 
$$ V_{g,n+2} \asymp \sqrt{n}\; V_{g-1,n+4}.$$
\end{itemize}
\end{subsection}
\begin{subsection}{}
It is crucial for the applications in $\ch \ref{Rs}$ to understand the behaviour of the polynomial $V_{g,1}(L)$ as $g \rightarrow \infty.$
We know that in general $ V_{g,1}(L) \leq e^{\frac{L}{2}} V_{g,1}. $
 In view of ($\ref{re2}$), the estimates we prove in Theorem $\ref{first}$ (a), imply that 
if $L \ll g$ 
$$ V_{g,1}(L) \asymp e^{\frac{L}{2}} V_{g,1}.$$

To simplify the notation, let 
$$[{\bf x}]_{g,n}:= [\tau_{x_{1}},\ldots, \tau_{x_{n}}]_{g,n},$$
where ${\bf x}=(x_{1},\ldots,x_{n}).$
Also, given ${\bf b}=(b_{1},\ldots,b_{l})$ and ${\bf c} =(c_{1},\ldots,c_{m}),$ let $${\bf b} \oplus {\bf c} =(b_{1}+c_{1}-1,b_{2},\ldots,b_{l},c_{2},\ldots, c_{m}). $$
%\begin{lemm}\label{estimate1}
%Given $k_{1},k_{2},$ and $g\geq 0,$ we have 
%$$ \sum _{g_{1}+g_{2}=g} [x_{1},\ldots, x_{k_{1}}]_{g_{1},k_{1}} [y_{1},\ldots, y_{k_{2}}]_{g_{1},k_{2}} \leq [ x_{1}+x_{2}-1, x_{2},\ldots,x_{k_{1}}, y_{2},\ldots, y_{k_{2}}]_{g-1,k_{1}+k_{2}-1} $$
%\end{lemm}
The following lemma plays an important role in the proof of Theorem \ref{11}.
\begin{lemm}\label{upperL}
In terms of the above notation, for ${\bf x}=(x_{1},\ldots, x_{l}),$ and ${\bf y}=(y_{1},\ldots, y_{m})$, we have 
\begin{equation}\label{usesimple}
\sum_{g_{1}+g_{2}+1=g \atop g_{2}+1\geq g_{1} \geq 0} [{\bf x}]_{g_{1},l} \times [{\bf y}]_{g_{2}+1,m} 
\leq C[ {\bf x} \oplus {\bf y}]_{g-1,n-1}, 
\end{equation}
where $n=l+m$, and $C>0$ is a constant independent of $g,$ $n$, ${\bf x}$, and ${\bf y}.$
\end{lemm}
Note that $[x_{1},\ldots,x_{n}]_{g,n} > 0$ if and only if $x_{1}+\ldots+x_{n} \leq 3g-3+n.$ We remark that if at least one term on the left hand side of $(\ref{usesimple})$ is non-zero 
then $x_{1}+y_{1}-1+x_{2}+\ldots+x_{l}+y_{2}+\ldots+y_{m} \leq 3(g-1)-3+(n-1)$ and hence $[ {\bf x} \oplus {\bf y}]_{g-1,n-1}>0.$\\

\noindent
{\bf Sketch of proof.} The proof is by induction on $2g-2+n.$
 The main idea is using $({\bf III})$ for ${\bf x}$ and ${\bf x} \oplus {\bf y}. $ 
 First, we can choose $C$ such that $\ref{usesimple}$ holds for $2g-2+n \leq 2.$
 
 Expand $[{\bf x} \oplus {\bf y}]_{g-1,n-1}$ and all the terms including $x_{i}'$s in $(\ref{usesimple})$ using the recursive relation $({\bf III}).$ 
 
 Roughly speaking, since all the terms in equation $({\bf III})$ and ($\ref{usesimple}$) are positive, it is enough to check that after 
 expanding both sides every term on the left hand side has a corresponding term on the right hand side.\\
Here we check this for the terms in $\mathcal{B}$ defined by ($\ref{B}$). By the definition, in $\mathcal{B}_{\bf x}$ for ${\bf x}=(x_{1},\ldots x_{l})$ the coefficient of a term 
$[k_{1},k_{2},x_{3},\ldots x_{l}]_{g'}$ with $g'=g_1-1$ and $k_{1}+k_{2}=L+x_{1}-2$ is equal to $a_{L}$. Now using the induction hypothesis for ${\bf x}'= (k_{1},k_{2},x_{3},\ldots x_{l})$ and ${\bf y}$ implies that the contribution of the terms corresponding to ${\bf x}'$ is 
$ \leq C a_{L} [{\bf x}' \oplus {\bf y}]_{g-2,n}.$ On the other hand, when we apply $({\bf III})$ to $[{\bf x} \oplus {\bf y}]_{g-1,n-1}$ the coefficient of the term 
$ [{\bf x}' \oplus {\bf y}]_{g-2,n} $ in $\mathcal{B}$ is equal to $a_{L'}$ where $L'= k_{1}+y_{1}-1+k_{2}- (x_{1}+y_{1}-1) +2= k_{1}+k_{2}-x_{1}+2=L.$ 
\hfill $\Box$\\

By applying this result to $m=\ell=1,$ ${\bf x}=(1),$ and ${\bf y}=(0)$, Lemma $\ref{obser}$(a) implies:

\begin{coro}
As $g \rightarrow \infty$
\begin{equation}\label{ssimple}
\sum _{i=1}^{g-1} V_{i,1} V_{g-i,1} =O(V_{g-1,1})=O(\frac{V_{g}}{g}).
\end{equation}
\end{coro}
\noindent
{\bf Remark.} Similarly, by induction on $2g+n$ and using ${\bf III}$ one can show that, for $d\geq 1$
$$ \sum_{a_{1}+a_{2}=d} [a_{1},a_{2},x_{1},\ldots, x_{n}] \leq [d-1,0,x_{1},\ldots,x_{n}] $$
and 
$$[m, x_{1},\ldots,x_{n-1}]_{g,n} \leq [m-1,0,0,x_{1},\ldots, x_{n-1}]_{g-1,n+2}.$$
We skip the proofs since we won't need these inequalities in this paper.\\
Now we can prove the main result of this section:
\begin{theo}\label{first}
Let $n \geq 0.$ 
\begin{itemize}
\item {\bf a):} For any $k \in {\mathbb N}$
 %\begin{equation}\label{this}
$$\frac{ [\tau_{k},\tau_{0},\ldots,\tau_{0}]_{g,n+1}} {V_{g,n+1}}=1+ O(1/g),$$ 
%\end{equation}
as $g \rightarrow \infty.$

\item{\bf b):} $$\frac{V_{g,n+1}} {2gV_{g,n}}= 4\pi^{2}+ O(1/g),$$
\item {\bf c):}
$$
 \frac{V_{g,n}} {V_{g-1,n+2}}=1+O(1/g).$$
 \end{itemize}
\end{theo}
\noindent
{\bf Remark.} 
\begin{itemize}
\item These estimates are consistent with the conjectures on the growth of Weil-Petersson volumes in \cite{Z:con}; we remark that the statements had been predicted by Peter Zograf.
\item Following the ideas used in the proof, one can show that 
$$\frac{V_{g,n+1}} {2gV_{g,n}}= 4\pi^{2}+ \frac{a_{1,n}}{g}+\ldots+\frac{a_{k,n}}{g^{k}}+O(\frac{1}{g^{k+1}}),$$
and
$$\frac{V_{g,n}} {V_{g-1,n+2}}=1+ \frac{b_{1,n}}{g}+\ldots+\frac{b_{k,n}}{g^{k}}+O(\frac{1}{g^{k+1}}).$$
However, in general it is not easy to calculate $a_{i,n}$ and $b_{i,n}'$s.
\end{itemize}
\noindent 
{\bf Proof of Theorem \ref{first}.} Fix $n \geq 0.$
Applying Lemma $\ref{upperL}$ for $(1,\ldots, 0)$ and $(0,\ldots,0)$ implies that as $g \rightarrow \infty$

\begin{equation}\label{bound1}
\sum_{g_{1}+g_{2}=g\; \atop \{1,\ldots,n\}=I _{1}\amalg I_{2}} V_{g_{1},|I_{1}|+1} \times V_{g_{2},|I_{2}|+1} = O(V_{g-1,n+1})
\end{equation}
\begin{equation}\label{bound2}
\sum_{g_{1}+g_{2}=g\; \atop \{1,\ldots,n\}=I_{1} \amalg I_{2}} V_{g_{1},|I_{1}|+2} \times V_{g_{2},|I_{2}|+2} = O(V_{g-1,n+3}).
\end{equation}

First, we assume that $n \geq 1.$
Then by the inequalities of $(\ref{yekido})$
$$ V_{g-1,n+2}=O(\frac{V_{g,n+1}}{g}), \; \; \mbox{and}\;\; V_{g,n}=O(\frac{V_{g,n+1}}{g}). $$
and from $(\ref{bound1})$ we get that 
$$ \sum_{g_{1}+g_{2}=g\; \atop \{1,\ldots,n\}=I _{1}\amalg I_{2}} V_{g_{1},|I_{1}|+1} \times V_{g_{2},|I_{2}|+1} = O(V_{g-1,n+1})=O(\frac{V_{g,n+1}}{g}). $$
Now by comparing the contributions of $\mathcal{A}_{\bf d}$, $\mathcal{B}_{\bf d},$ and $\mathcal{C}_{\bf d}$ for ${\bf d}=(k,0,\ldots,0)$ and $(0,\ldots,0)$ in $({\bf III})$, 
and Lemma $\ref{aa}$ we get
\begin{equation}\label{1g}
|\frac{ [\tau_{k},\tau_{0},\ldots,\tau_{0}]_{g,n+1}} {V_{g,n+1}}-1| \leq c_{0} \frac{k^2}{g},
\end{equation}
where $c_{0}$ is a universal constant independent of $g$ and $k$.

We use the following elementary observation to prove $({\bf b})$ for $n\geq 1$: \\
\noindent
{\bf Elementary fact}.
{\it Let $\{r_{i}\}_{i=1}^{\infty}$ be a sequence of real numbers and $\{k_{g}\}_{g=1}^{\infty}$ be an increasing sequence of positive integers.
Assume that for $g\geq 1$, and $i \in {\mathbb N}$, $0 \leq c_{g,i} \leq c_{i},$ and 
 $\lim_{g\rightarrow \infty} c_{g,i}=c_{i}.$
If $\sum_{i=1}^{\infty} | c_{i} r_{i}| < \infty,$ then 
\begin{equation}\label{fact}
\lim_{g\rightarrow \infty} \sum_{i=1}^{k_{g}} r_{i} c_{g,i} = \sum_{i=1}^{\infty} r_{i} c_{i}.
\end{equation}}

Now, let 
$$r_{i}= (-1)^{i} \frac{\pi^{2i} (i+1)}{(2i+3)!}, \; k_{g}=3g-3+n\;, c_{i}=1\;\; \mbox{and}\;\; c_{g,i}= \frac{[\tau_{i+1}\tau_{0}\ldots \tau_{0}]_{g,n}}{V_{g,n+1}}.$$
By $(\ref{fact}),$ and $({\bf II})$ for ${\bf d}=0$ we get 
$$
\lim _{ g \rightarrow \infty} \frac{2 (2g-2+n) V_{g,n}}{V_{g,n+1}}= \frac{1}{3!} -\frac{2 \pi^{2}}{5!}+\ldots+ (-1)^{L} (L+1) \frac{\pi^{2L}}{(2L+3)!}+\ldots = \frac{1}{2\pi^{2}}.$$
In fact, similarly $(\ref{1g})$ implies that 
$$\frac{2 (2g-2+n) V_{g,n}}{V_{g,n+1}}=\frac{1}{2\pi^2}+O(\frac{1}{g}).$$
On the other hand, from $({\bf I})$ and $(\ref{bound2})$ we get that for $n \geq 2:$
$$\lim_{g\rightarrow \infty} \frac{V_{g,n}} {V_{g-1,n+2}}= 1+O(1/g).$$
Now it is easy to check that 
$$V_{g,1}=\frac{1}{g} V_{g,2}( \frac{1}{4\pi^2} (1-O(1/g)) \;, \; V_{g-1,3}=\frac{1}{g}V_{g-1,4}( \frac{1}{4\pi^2} (1-O(1/g))$$
and $$V_{g,2}=V_{g-1,4}(1+O(1/g))$$ imply 
$$\frac{V_{g,1}}{V_{g-1,3}}=1+O(1/g).$$
In other words, $({\bf b})$ for $n=1$ and $n=2$ proves $({\bf c})$ for $n=1$.

We remark that $(\ref{zero})$ implies $({\bf b})$ for $n=0.$
 Finally $({\bf b})$ for $n=0$ and $n=1$ implies $({\bf c})$ for $n=0$. 
\hfill $\Box$\\
It is easy to check that if $\{k_{n}\},$ is a bounded sequence $|k_{n}|< c$ 
$$\frac{1}{m} g^{-c} < \prod_{n=1}^{g} (1+\frac{k_{n}}{n}) <m g^{c},$$
where $m$ is independent of $g.$
Hence, Theorem $\ref{first}$ implies the following 
\begin{coro}\label{power}
There exists $M>0$ such that: 
\begin{equation}\label{claimm1}
g^{-M} \cC_{g,n}< V_{g,n}< g^{M} \cC_{g,n}, 
\end{equation}
where
$$ \cC_{g,n}=(4\pi^2)^{2g+n-3} (2g-3+n)! \frac{1}{\sqrt{g \pi}}.$$
\end{coro}
Finally, $(\ref{claimm1})$ implies that 
\begin{equation}\label{uuse}
\sum _{i=r+1}^{g/2} V_{i,1} \times V_{g-i,1} \asymp \frac{V_{g}}{g^{2r+1}}.
\end{equation}
%\asymp V_{g,1}/g^{2} 
 A simple calculation shows that 
 \begin{equation}\label{simple}
 \sum _{g_{0}+g_{1}=g+1-k, \atop r\leq g_{0} \leq g_{1} } e^{C g_{0}} g_{0} \frac{\cC_{g_{0},k} \cC_{g_{1},k}}{\cC_{g}}=O(\frac{1}{g^{r}}), 
 \end{equation}
where $C=2 \ln(2).$
Therefore, we get the following estimate which will be used in the next section: 
\begin{coro}\label{b0}
Let $k \geq 0,$ $1>\beta>0$, and $c>0.$ Then as $g \rightarrow \infty$
$$\sum _{g_{0}+g_{1}=g+1-k } e^{C g_{0}+ c g_{0}^{\beta}} g_{0} V_{g_{0},k} V_{g_{1},k}= O (\frac{V_{g}}{g}), $$
where $C=2 \ln(2).$
\end{coro}
\end{subsection}
 \end{section}
\begin{section}{Random Riemann surfaces of high genus}\label{Rs}
In this section, we apply the asymptotic estimates on the volume polynomials 
to study the geometric properties of random hyperbolic surfaces; in particular, we are interested in the length of the shortest simple closed geodesic of a given 
combinatorial type, diameter and the Cheeger constant of a random surface. See \cite{B:M} for more in the case of random hyperbolic surfaces constructed by random 
trivalent graphs.
\begin{subsection}{}\label{NR}
{\bf Notation.} Recall that the mapping class group $\operatorname{Mod}_{g,n}$ acts naturally
on the set $\mathcal{S}_{g,n}$ of isotopy classes of simple closed curves on $S_{g,n}$: Two simple closed curves $\alpha_{1}$ and 
$\alpha_{2}$ are of the same {\it type} if and only
if there exists $g \in \operatorname{Mod}_{g,n}$ such that 
 $g \cdot \alpha_{1}=\alpha_{2}.$ The type of a simple closed
curve is determined by the topology of $S_{g,n}-\alpha$, the surface that we get
by cutting $S_{g,n}$ along $\alpha$.\\
Let 
 $$\mathcal{S}_{k}^{m}=\{ \gamma=\gamma_{1}+\ldots+\gamma_{k} |\; \gamma_{i} \in \mathcal{S}_{g,n}\; \mbox{distinct} \; S-\gamma= S_{1} \cup S_{2}, |\chi(S_{1})|=m\}.$$
 Note that each $\mathcal{S}_{k}^{m}$ is invariant under the action of the mapping class group.
To simplify the notation, let $\widetilde{\gamma_{0}}$ be a non-separating simple closed curve on $S_{g},$ and $\widetilde{\gamma_{i}}$ be a separating simple closed curve
on $S_{g}$ such that $$S_{g}-\widetilde{\gamma_{i}}=S_{i,1}\cup S_{g-i,1}.$$ That is $\tilde{\gamma_{i}} \in \mathcal{S}_{1}^{2i-1},$ and $\widetilde{\gamma_{0}}= \mathcal{S}_{1}^{2g-2}.$

Consider the counting function 
$$N_{\alpha}(\cdot, \cdot) : {\mathbb R}_{+} \times \mathcal{M}_{g} \rightarrow {\mathbb R}_{+}$$
defined by 
$$N_{\alpha} (L,X)= |\{\gamma | \gamma \in \alpha\cdot \Mod_{g}, \ell_{\gamma}(X) \leq L\} |.$$
Let
\begin{equation}\label{def}
F^{L}_{i}(X)=| \{\gamma | \gamma \in \mathcal{O}_{\gamma_{i}}, \ell_{\gamma}(X) \leq L\} |.
\end{equation}
Using Theorem $\ref{integrate}$ and the estimates proved in \ch $\ref{Asym}$, we can show that if $L$ is fixed
$$ \int _{\mathcal{M}_{g}} F_{0}^{L}(X) \; dX= \int_{0}^{L} t V_{g-1,2}(t,t) dt \asymp e^{L}L^{2} V_{g},$$
but 
$$ \int _{\mathcal{M}_{g}} F_{1}^{L}(X) \; dX \asymp \frac{e^{L/2}L^{3} V_{g}}{g},$$
as $g \rightarrow \infty.$
Similar estimates hold when $L$ is {\it much smaller} than $g.$ We remark that since the 
number of closed geodesics of length $\leq L$ on 
$X \in \mathcal{M}_{g}$ is at most $e^{L+6} (g-1)$ (see Lemma $6.6.4$ in \cite{B:book}), we can not
expect the similar bounds to hold in general. However, in general we have 
\begin{equation}\label{vs}
\int _{\mathcal{M}_{g}} F_{0}^{L}(X) \; dX= O(e^{L}L^{2} V_{g}) \; ,\mbox{and}  \int _{\mathcal{M}_{g}} F_{1}^{L}(X) \; dX =O( \frac{e^{L/2}L^{3} V_{g}}{g}).
\end{equation} 

\end{subsection}
%In this section, we will discuss some geometric properties of such a surface through discussing
%the behavior of the length of the shortest simple geodesic of a given combinatorial type on $X.$
%A REMARK HERE :1) MAXIMUM VOLUME IS FOR THE CONNECTED GUYS! 2) DISCTONNECTED GUYS ARE O(EPSILON2/G) VA AXIMUM TEDAAD 
%G SO THEY DON'T CONTRIBUTE TO THE COUNTING OF SHORT CURVES. 1/E WILL HAVE A SHORT CURVE---AT LEAST THIS-----\\ 
\begin{subsection}{}\label{sh}
\noindent 
{\bf Systoles and injectivity radius.}
Let $$\mathcal{M}_{g,n}^{\epsilon}=\{ X\;| \; \exists \gamma, \ell_{\gamma}(X) \leq \epsilon\}\subset \mathcal{M}_{g,n}.$$
The set $\mathcal{M}_{g,n}-\mathcal{M}_{g,n}^{\epsilon}$ of hyperbolic surfaces with lengths of closed geodesics bounded below by a constant $\epsilon > 0$ is 
a compact subset of the moduli space $\mathcal{M}_{g,n}$.
\begin{theo}\label{thin}
Let $n\geq 0.$ There exists $\epsilon_{0}>0$ such that for any $\epsilon< \epsilon_{0}$ 
$$\vo_{wp}(\mathcal{M}_{g,n}^{\epsilon}) \asymp \epsilon^2 \vo_{wp}(\mathcal{M}_{g,n})$$
as $g \rightarrow \infty.$
\end{theo}
\noindent
{\bf Proof.}
Here we sketch the proof for the case of $n=0.$
Fix $\epsilon$ such that no two simple closed geodesics of length $\leq \epsilon$ could meet.
Consider the function
$$F^{\epsilon}(X)=N (\epsilon,X) = F^{\epsilon}_{0}(X)+\ldots F^{\epsilon}_{g/2}(X),$$
as defined in $(\ref{def}).$
Then by Theorem \ref{integrate}, we have
$$ \vo_{wp}(\mathcal{M}_{g}^{\epsilon}) \leq \int _{\mathcal{M}_{g}} F^{\epsilon}(X) \; dX \leq$$
$$\leq \sum _{i=1}^{g/2} \int_{0}^{\epsilon} t \vo_{wp}(\mathcal{M} (S_{g}-\gamma_{i}, t,t)) \;dt \;+ \int_{0} ^{\epsilon} t \vo_{wp}(\mathcal{M}_{g-1,2}(t,t))\; dt $$
On the other hand, by $(\ref{upper})$ we know that if $t$ is small enough for $i \geq 1,$
$$ \vo_{wp}(\mathcal{M}(S_{g}-\gamma_{i}, t,t))\leq 2 V_{i,1} \times V_{g-i,1},$$
and 
$$ \vo_{wp}(\mathcal{M}_{g-1,2}(t,t))\leq 2 V_{g-1,2}. $$
Hence, when $\epsilon$ is small (independent of $g$), from $(\ref{uuse})$ and $(\ref{ssimple})$ we get 
$$ \vo_{wp}(\mathcal{M}_{g}^{\epsilon})=O (\epsilon^2 (\sum_{i=1}^{g/2} V_{i,1} V_{g-i,1}+V_{g-1,2}))= O (\epsilon^2 V_{g}). $$
Next, we prove that the volume of the locus with a {\it non-separating} short simple closed geodesic of length $\leq \epsilon$ is asymptotically positive. 

Since 
$ \int _{\mathcal{M}_{g}} F_{0}^{\epsilon}(X) \; dX \asymp \epsilon^2 V_{g}$ in order to get a lower bound $\vo_{wp}(\mathcal{M}_{g,n}^{\epsilon})$ we need to prove an upper bound for the volume of the locus where $F_{0}^{\epsilon}(X) \geq k$ for $k \geq 2.$ 
In fact 
$$ \int _{\mathcal{M}_{g}} F_{0}^{\epsilon}(X)= \sum_{k=1}^{\infty} \vo_{wp}(\{X | F_{0}^{\epsilon} \geq k \}). $$
Note that by Lemma $\ref{obser}$
$$V_{g-1,n+4} \leq V_{g,n+2}, $$
and if $\sum_{i=1}^{s} (2g_{i}-2+k_{i})=2g-2$ with $s \geq 2,$ and $k_{i} \geq 2$ for $1 \leq i\leq s $
 
\begin{equation}\label{other}
\prod V_{g_{i},k_{i}}= O (\frac{V_{g}}{g^{2}}). 
\end{equation}
Let 
$$\mathcal{U}=\{ X \; | \exists \gamma_{1},\ldots \gamma_{l} \in \Mod \cdot \gamma_{0}, \ell_{\gamma_{i}}(X) \leq \epsilon, \; S-\cup \gamma_{i}\; \mbox{is disconnected}\; \} \subset \mathcal{M}_{g}.$$
Since $F_{0}^{\epsilon}(X) \leq 3g-3$, from $(\ref{other})$ 
$$ \int _{\mathcal{U}} F_{0}^{\epsilon}(X) \; dX \leq O (\frac{\epsilon^{4} V_{g}}{g}).$$
On the other hand, by using the same argument for 
$$F_{0,k}^{\epsilon}(X)=| \{ \{\gamma_{1},\ldots, \gamma_{k}\} | \gamma_{i} \; \mbox{non-separating} \; \ell_{\gamma_{i}}(X) \leq \epsilon\} |$$ and applying Theorem \ref{integrate}, 
we get 
$$ \vo_{wp}(\{X | F_{0, k}^{\epsilon} \geq 1\}- \mathcal{U})= \vo_{wp}(X | F^{\epsilon} \geq k \} - \mathcal{U}) \leq c \frac{\epsilon^{2k}e^{\epsilon k}}{k!},$$
where $c$ is a constant independent of $g$ and $k$.
Therefore if $\epsilon>0$ is small enough 
$$\sum_{k=2}^{\infty} \vo_{wp}(\{X | F_{0}^{\epsilon}(X) \geq k\}) \leq \epsilon^4 \vo_{wp}(\mathcal{M}_{g})$$
which implies the result.
\hfill $\Box$\\
Let $$f(X)= \sum_{\ell_{\alpha}(X) \leq 1} \frac{1}{\ell_{\alpha}(X)}.$$
Then using Theorem $\ref{integrate}$
$$\int_{\mathcal{M}_{g}} f(X) \;dX= \int_{0}^{1} V_{g-1,2}(t,t) dt+ \sum_{i=1}^{g/2} \int_{0}^{1} V_{g-i,1}(t) V_{i,1}(t) dt \asymp V_{g}$$
and hence Theorem $\ref{thin}$ implies that :
\begin{coro}\label{as}
As $g \rightarrow \infty$
$$
\int_{\mathcal{M}_{g}} \frac{1}{\ell_{sys}(X)} \; dX \asymp V_{g}. $$
\end{coro}
\end{subsection}
\begin{subsection}{}\label{ss}
\noindent
{\bf Behavior of separating simple closed geodesics.}
By \cite{S:S} there exists a positive constant $C>0$ such that 
every closed surface $X$ of genus $g\geq 2$, $\ell_{sys}^{s}(X) \leq C \log(g).$ 
We show that as $g \rightarrow \infty$ $\ell_{sys}^{s}(X)$ is generically at least of $(2-\epsilon)\log(g)$.
Moreover generically if a separating curve $\gamma$ satisfies $\ell_{\gamma}(X) < C \log(g)$ then $S_{g}-\gamma= S_{g_{1}} \cup S_{g_{2}}$ with $g_{1}=O(1).$ 
\begin{theo}
Let $0<m < 2$ then
$$ \PP^g_{wp} (\ell^{s}_{sys}(X)< m \log(g))=O(\log(g) g^{(m/2-1)})),$$
and
$$\EE^g_{X\sim wp} (\ell^{s}_{sys}(X)) \asymp \log(g)$$
as $g \rightarrow \infty.$
\end{theo}
\noindent
{\bf Proof.} 
Note that 
\begin{equation}\label{mvs}
\PP^g_{wp} (\ell^{s}_{sys}(X)< L ) \leq \frac{e^{L/2}L^{3}}{g}+ \sum _{i=2}^{g/2} \frac{ e^{L} V_{i,1}\times V_{g-i,1}}{V_{g}}. 
\end{equation}
On the other hand, by $(\ref{uuse})$
$$ \sum _{i=2}^{g/2} e^{L} V_{i,1}\times V_{g-i,1} =O(\frac{V_{g}}{g^{3}})$$
which implies the result. 
\hfill $\Box$\\
\end{subsection}
\begin{subsection}{}\label{Ie}
\noindent
{\bf Injectivity radius and embedded balls.}
Let $\operatorname{Inj}(x)$ denote the injectivity radius of $x \in X$. We show that on a generic $X \in \mathcal{M}_{g}$ almost every point $x \in X$ has $\operatorname{Inj}(x) \geq \frac{1}{6} \log(g).$ 
By the definition of the injectivity radius, corresponding to each $x,$ there exists a simple closed curve 
$\gamma_{x}$ of length $\leq 2 \operatorname{Inj}(x)$ such that the distance of $x$ from the geodesic representative of $\gamma_{x}$ is at most $2 \operatorname{Inj}(x).$ 
Also, let $N(L,X)$ be the number of simple closed geodesics of length $\leq L$ on $X$. Then 
\begin{itemize}
\item It is easy to see from $(\ref{vs})$ and $(\ref{mvs})$ that  as $g \rightarrow \infty $ 
\begin{equation}\label{dd}
\PP^g (\{X\; | \; N (\log(g)/3,X) \geq g^{1/3+1/4}) =O(g^{-1/4}). 
\end{equation}
\item A simple calculation shows that given a simple closed geodesic $\gamma$ of length $\leq \log(g)/3$ the volume of the locus on $X$ with $\gamma_{x}=\gamma$ is at most $g^{1/3} \log(g).$
\end{itemize}
Therefore, for a generic point in $X \in \mathcal{M}_{g}$ (defined by $(\ref{dd})$)
$$ \vo(\{ x \in X\; | \operatorname{Inj}(x) \leq \frac{1}{6}\log (g)\})=O(g^{11/12} \log(g)).$$
Hence, we have: 
\begin{theo}
As $g \rightarrow \infty$
$$\PP^g_{wp}(\operatorname{Emb}(X) \leq C_{E}\log(g)) \rightarrow 0,$$
$$\EE^g_{X\sim wp} (\operatorname{Emb}(X)) \asymp \log(g), $$
where $C_{E}=\frac{1}{3}.$ 
\end{theo}
\end{subsection}
\begin{subsection}{}\label{ch}
\noindent
{\bf Cheeger constants and isoperimetric inequalities.}
 Recall that the Cheeger constant of $X$ is defined by 
$$h(X)= \inf \frac{\ell(A)}{ \min \{\operatorname{Area}(X_{1}), \operatorname{Area}(X_{2})\}}$$
where the infimum is taken over all smooth $1$-dimensional submanifolds of $X$ which divide it into two disjoint submanifolds $X_{1}$ and $X_{2}$ such that 
$X-A=X_{1} \cup X_{2}$ and $A \subset \partial(X_{1}) \cap \partial(X_{2}).$\\
\noindent
We remark that:
\begin{itemize}
\item In fact, by an observation due to Yau, we may restrict $A$ to a family of curves for which $X_{1}$ and $X_{2}$ are connected. See \cite{B:book}.
\item By a result of Cheng \cite{C}, $$h(X) \leq 1+ \frac{16\pi^2}{\operatorname{diam}(X)}.$$
Therefore, there is an upper bound for the Cheeger constant which tends to $1$ as $g(X) \rightarrow \infty.$
See also \ch III and \ch X in \cite{C:book}.
\end{itemize}

Given $i \leq g$
$$H_{i}(X)= \inf \frac{\ell_{\alpha}(X)}{ \min \{\operatorname{Area}(X_{1}), \operatorname{Area}(X_{2})\}}$$
where $\alpha= \cup_{j=1}^{s} \alpha_{j}$ is a union of simple closed {\it geodesics} on $X$ with $X-\alpha= X_{1} \cup X_{2},$ and $X_{1}$ and $X_{2}$ are connected
subsurfaces of $X$ such that $| \chi(X_{1})|=i \leq | \chi(X_{2})|$. Here $| \chi(X_{1})|= 2g_{1}-2+s.$
 Now we define the {\it geodesic} Cheeger constant of $X$ by
 $$H_{i}(X)= \min_{i \leq g} H_{i}(X).$$
In general, by the definition $$H(X) \geq h(X),$$ but the inequality is not sharp. However, using the following basic properties of perimeter minimizers allows us to get a lower bound for 
$h(X)$ in terms of $H(X)$.

Recall that in a compact hyperbolic surface, there exists a perimeter minimizer among regions of prescribed area bounded 
by embedded rectifiable curves; it consists of curves of equal constant curvature. Moreover, by a result of Adams and Morgan \cite{AM:I}:
 \begin{theo}
 For given area $0 < A< 4 \pi g $, a perimeter-minimizing system of embedded rectifiable curves bounding a region $R$ of 
area $A$ consists of a set of curves of one of the following four types :

\begin{enumerate}
\item a circle, 
\item horocycles around cusps, 
\item two Òneighboring curvesÓ at constant distance from a geodesic, bounding an 
annulus or complement, 
\item geodesics or single Òneighboring curves.Ó 
\end{enumerate}
All curves in the set have the same constant curvature.

\end{theo}
In fact in the case of a circle or neighboring curves $h$ is strictly bigger than $1$. 
On the other hand, by a simple calculation (see Lemma 2.3 \cite{AM:I}) if a neighboring curve of length $L$ and curvature $\kappa$ at distance $s$ from a geodesic of length $\ell$, enclosing area $A$, then
 $$A =\ell \sinh(s),\;\;\;L = \ell \cosh(s) \;\;, \mbox{and} \;\; \kappa = \tanh (s).$$

Therefore by using basic isoperimetric inequalities for hyperbolic surfaces, we have: 
\begin{prop}\label{hH}
Let $X \in \mathcal{M}_{g}$ be a hyperbolic surface of genus $g$.
Then 
$$ h(X) \geq \frac{H(X)}{H(X)+1}$$
\end{prop}
Now, we can show: 
\begin{theo}\label{Cheeger}
As $g \rightarrow \infty$ 
$$\PP^g_{wp}\left (h(X) \leq \frac{\ln(2)}{\pi+\ln(2)}\right) \rightarrow 0,$$
and
\begin{equation}\label{ah}
\int_{\mathcal{M}_{g}} \frac{1}{h(X)} \; dX \asymp V_{g}.
\end{equation}
\end{theo}

Let 
$$ \mathcal{W}_{k}^{2m-1}(L)= \vo _{wp} (\{ X\in \mathcal{M}_{g}\; | \exists \gamma \in \mathcal{S}_{k}^{2m-1} \;\; ,\ell_{\gamma}(X) \leq L\}.$$ 
Let $\gamma_{k}^{2m-1}=\gamma_{1}+\ldots \gamma_{k} \in \mathcal{S}_{k}^{2m-1}$ (see \ch $\ref{NR}$). By Theorem $\ref{integrate}$ for $N_{\gamma_{k}^{2m-1}}(X,L)$ we have:
$$\mathcal{W}_{k}^{2m-1}(L) \leq e^{L} V_{m,1} \times V_{g-m,1} \times \int_{L_{1}+\ldots L_{k}\leq L} \frac{1}{k!} L_{1}\cdots L_{k}\; dL_{1}\cdots dL_{k},$$
and 
$$\mathcal{W}^{2m-1}(L) \leq e^{L} V_{m,1} \times V_{g-m,1} \sum_{k=1}^{2m} \int_{L_{1}+\ldots L_{k}\leq L} \frac{1}{k!} L_{1}\cdots L_{k}\; dL_{1}\cdots dL_{k}.$$
On the other hand,
since 
$$ \int_{L_{1}+\ldots L_{s}\leq L} L_{1}\cdots L_{s}\; dL_{1}\cdots dL_{s}= \frac{L^{2s}}{(2s)!},$$ 
and 
$$ \sum_{s=1}^{\infty} \frac{L^{2s}}{s! (2s)!}= O( e^{ \frac{3}{2} L^{2/3}})$$
we get 
\begin{equation}\label{os}
\mathcal{W}^{2m-1}(L) \leq e^{L+\frac{3}{2}L^{2/3}} V_{m,1} \times V_{g-m,1}.
\end{equation}
As before, let 
$$H_{k}(X)= \inf \frac{\ell_{\alpha}(X)}{\pi k }$$
where $\alpha= \cup_{i} \alpha_{i}$ is a union of simple closed {\it geodesics} on $X$ with $X-\alpha= X_{1} \cup X_{2},$ and $X_{1}$ and $X_{2}$ are connected
subsurfaces of $X$ such that $|\chi (X_{1})| = k < | \chi(X_{2})|. $
Recall that by Lemma \ref{obser}, for $n\geq 0$ 
$$V_{g-1,n+4} \leq V_{g,n+2}.$$
Hence, from $(\ref{os})$ we get: 
\begin{lemm}\label{b1}
Let $m=2m_{1}-2+n_{1} \leq 2g-2, $
where $n_{1} \in \{0,1\}.$
Then 
$$\vol_{wp} ( \{ X | H_{m}(X) \leq C \} = O(m e^{\pi\cdot m \cdot C+c m^{2/3}} V_{m_{1}+n_{1}-1,2-n_{1}} V_{2g-2-m+n_{1},2-n_{1}}), $$
where $c$ is a constant independent of $g$.
\end{lemm}
\noindent
{\bf Proof of Theorem \ref{Cheeger}.} 
Lemma $\ref{b1}$ and Corollary $\ref{b0}$ imply that as $g \rightarrow \infty$
$$\PP^g_{wp}\left (H(X) \leq \frac{\ln(2)}{\pi}\right) \rightarrow 0.$$
 Therefore, in view of Proposition $\ref{hH}$, we get the result. 
 Corollary $\ref{as}$ implies the second part of the theorem.
 \hfill $\Box$\\
 Moreover, we have:
 \begin{theo}
Let $s_{g}$ be a sequence such that $\lim_{g \rightarrow \infty} \frac{s_{g}}{g}=0.$
Given $M>0,$ 
$$\PP^g_{wp} (H_{s_{g}}(X) \leq M) \rightarrow 0,$$
as $g \rightarrow \infty.$
\end{theo}
\end{subsection}

\begin{subsection}{}\label{di}
\noindent
{\bf Diameter.}
It is known that the diameter of a Riemannian manifold of constant curvature $-1$ satisfies:
\begin{equation}\label{relation}
\operatorname{diam}(X) \leq 2 (r_{0}+ \frac{1}{h}\log (\frac{\operatorname{Vol}(X)}{2 B(r_{0})})),
\end{equation}
where $r_{0}>0$ and $B(r_{o})$ is the infimum of the volume of a ball of radius $r_{0}$ in $X$.
%(ADD REF.)????????????????????????????????????????????????????????????????
Using this result, we get:
\begin{theo}\label{dia}
As $g \rightarrow \infty$ 
$$\PP^g_{wp}(\operatorname{diam}(X) \geq C_{d}\log(g)) \rightarrow 0,$$
and
$$\EE^g_{X\sim wp} (\sqrt{\operatorname{diam}(X)}) \asymp \sqrt{\log(g)},$$
where $C_{d}=5.$
\end{theo}
\noindent
{\bf Proof.} From the proof of Theorem $\ref{thin}$ 
$$ \PP^g (\{ X \; | \; \exists \gamma \; \ell_{\gamma}(X) \leq \frac{1}{\log(g)}\}) = O (\frac{1}{\log(g)^2}).$$
Therefore the first part of this theorem is a direct consequence of $(\ref{relation})$ and Theorem $\ref{Cheeger}.$
Next we need to prove that 
$$\EE^g_{X\sim wp} (\sqrt{\operatorname{diam}(X)})=O(\sqrt{\log(g)}).$$
Since (as in Corollary $\ref{as}$) 
$$ \int_{\mathcal{M}_{g}}| \log(\ell_{sys}(X))| dX \asymp V_{g}$$
we have 
$$\int_{\mathcal{M}_{g}} \sqrt{\frac{|\log(\ell_{sys}(X))|}{h(X)}} \;dX \asymp V_{g}.$$

%In order to prove this claim, note that by Proposition \ref{hH} when $h(X)$ is very small
%$h(X) \asymp H(X).$ 
Now the second part follows from Corollary $\ref{as}$ and $(\ref{ah}).$

\hfill $\Box$\\
\end{subsection}
\end{section}

\bigskip
Department of Mathematics, Stanford University, Stanford CA 94305 USA;\\
mmirzakh@math.stanford.edu


\begin{thebibliography}{a}

 \bibitem[AC]{AC}
E.~Arbarello and M.~Cornalba.
 {\em Combinatorial and algebro-geometric cohomology
classes on the Moduli Spaces of Curves,}
J. Algebraic Geometry {\bf 5} (1996), 705--709.
 \bibitem[AM]{AM:I}
 C.~Adams and F. ~Morgan. {\em Isoperimetric curves on hyperbolic surfaces,} 
Proc. Amer. Math. Soc. {\bf 127} (1999), 1347-1356.

\bibitem [BS]{BS:S}
F. ~Balacheff and S. Sabourau
{\em Diastolic inequalities and isoperimetric inequalities on surfaces,}
 Annales Scientifiques de l'ƒcole Normale SupŽrieure, to appear

%\bibitem[Bo]{Bow:iden}
%B. ~Bowditch.
%{\em Markoff triples and quasifuchsian groups,}
 %Proceedings of the London Mathematical Society, Vol. 77 (1998), 697--736. 

 \bibitem[BM]{B:M}
 R.~Brooks and E.~Makover.
 {\em Random Construction of Riemann Surfaces,}
J. Differential Geom. {\bf 68} (2004), 121--157.

 \bibitem[Bu]{B:book}
P.~Buser. Geometry and spectra of compact Riemann surfaces, Birkh¬auser Boston, 1992. 
 
 \bibitem[BP]{BP:sys} 
 P.~Buser and P.~Sarnak.
 {\em On the period matrix of a Riemann surface of large genus,}
 Invent. Math. {\bf 117:1} (1994), 27--56.

 \bibitem[CP]{C:P}
W. Cavendish, H. Parlier.
{\em Growth of the Weil-Petersson Diameter of Moduli Space,} Preprint.

 
 \bibitem[Ch]{C:book}
 I.~Chavel, Eigenvalues in Riemannian Geometry, Academic Press, 1984.
 
 \bibitem[C]{C}
S. Cheng, { \em Eigenvalue comparison theorems and its geometric applications,}
 Math. Z. {\bf 143}, 289--297.


 \bibitem[DN]{DN:cone}
 N. Do and P. Norbury. {\em Weil-Petersson volumes and cone surfaces,} 
 Geom. Dedicata {\bf 141} (2009), 93--107.
 
 \bibitem[E]{E:M}
B. ~Eynard.
{\em Recursion between Mumford volumes of moduli spaces,}
Preprint.
 
 \bibitem[EO]{EO:WK}
B. ~Eynard and N. ~Orantin.
 {\em Invariants of algebraic curves and topological expansion,}
 Commun. Number Theory Phys. 1:2 (2007), 347--452. 
 \bibitem[Ga]{Ga}
A.~Gamburd.
{\em Poisson-Dirichlet distribution for random Belyi surfaces,}
Ann. Probab. {\bf 34:5} (2006), 1827--1848.

\bibitem[Go]{G:S}
W.~Goldman.
{\em The symplectic nature of fundamental groups of surfaces,}
Adv. Math. {\bf 54} (1984), 200--225.

\bibitem[Gr]{Gr:vol}
S.~Grushevsky.
 {\em An explicit upper bound for Weil-Petersson volumes of the moduli spaces of punctured Riemann surfaces,}
Mathematische Annalen. {\bf 321} (2001) 1, 1--13.

\bibitem[HM]{Harris:book}
 J.~ Harris and I.~ Morrison.
Moduli of Curves. Graduate Texts in Mathematics, vol 187, Springer-Verlag, 1998.
\bibitem[IT]{IT:book}
Y.~Imayoshi and M.~Taniguchi.
{\em An introduction to Teichm\"uller spaces,} Springer-Verlag, 1992.

\bibitem[Hu]{Hu} 
 Z. Huang, {\em On asymptotic Weil-Petersson geometry of Teichm\"uller space of Riemann surfaces.}
 Asian J. Math, {\bf 11} (2007), 459-484. 

\bibitem[IZ]{I:Z}
 C.~Itzykson and J.~Zuber. 
 {\em Combinatorics of the modular group. II. The Kontsevich integrals,} 
Internat. J. Modern Phys. A. {\bf 7} (1992), 5661--5705.

 \bibitem[KMZ]{KMZ:w}
 R.~Kaufmann, Y.~Manin, and D.~Zagier.
 {\em Higher Weil-Petersson volumes of moduli spaces of 
stable n-pointed curves,}
 Comm. Math. Phys. {\bf 181} (1996), 736--787.
 \bibitem[KL]{KL:W}
M. E.~Kazarian and S. K.~ Lando. 
{\em An algebro-geometric proof of Witten's conjecture,}
J. Amer. Math. Soc. {\bf 20} (2007), 1079--1089. 
 

 \bibitem[Ko]{Ko:int}
 M.~Kontsevich.
 {\em Intersection on the moduli space of curves and the matrix Airy function,}
Comm. Math. Phys. {\bf 147} (1992), 1-23.

 
 \bibitem[LX1]{LX:higher}
 K.~Liu and H.~Xu.
 {\em Recursion formulae of higher Weil-Petersson volumes}
 Int. Math. Res. Not. IMRN {\bf 5} (2009), 835--859.

 \bibitem[LX2]{LX:M}
 K. Liu, and H. Xu.
 {\em Mirzakharni's recursion formula is equivalent to the Witten-Kontsevich theorem,}
 Preprint.
 
\bibitem[MZ]{MZ}
Yu.~Manin and P. ~Zograf. {\em Invertible cohomological field theories and Weil-Petersson 
volumes,} Ann. Inst. Fourier {\bf 50:2} (2000), 519--535. 

\bibitem[MM]{M:M}
E.~Makover and J.~McGowan.
{\em The length of closed geodesics on random Riemann Surfaces,}
Preprint.

\bibitem[Mc]{M:S}
 G.~McShane. {\em Simple geodesics and a series constant over Teichm\"uller space.}
 Invent. Math. {\bf 132} (1998), 607--632.

 \bibitem[M1]{M:JAMS}
M.~Mirzakhani. {\em Weil-Petersson volumes and intersection theory on the moduli space of curves,}
 J. Amer. Math. Soc. {\bf 20:1} (2007), 1--23. 

\bibitem[M2]{M:In}
M.~Mirzakhani. {\em Simple geodesics and Weil-Petersson volumes of moduli spaces of bordered Riemann surfaces,}
 Invent. Math. {\bf 167} (2007), 179--222.

\bibitem[M3]{M:AB}
M.~Mirzakhani.
Random hyperbolic surfaces and measured laminations. In the tradition of Ahlfors-Bers. IV, 179--198, 
Contemp. Math., 432, Amer. Math. Soc., Providence, RI, 2007.

\bibitem[M4]{M:ICM}
M.~Mirzakhani.
On Weil-Petersson volumes and geometry of random hyperbolic surfaces. {\it Proceedings of ICM 2010.}


\bibitem[MS]{MuS}
Y.~Mulase and P.~Safnuk. {\em Mirzakhani's recursion relations, Virasoro constraints and the KdV hierarchy,}
Indian Journal of Mathematics {\bf 50} (2008), 189--228. 


\bibitem[OP]{OP}
 A.~Okounkov and R.~Pandharipande. {\em Gromov-Witten theory, Hurwitz theory, and matrix 
models, I,} Preprint. 
 \bibitem[Pe]{P:vol}
 R.~Penner.{\em Weil-Petersson volumes,}
 J. Differential Geom. {\bf 35} (1992), 559--608. 

 
 \bibitem[S1]{S:S1}
P.~ Schmutz.
{\em Geometry of Riemann surfaces based on closed geodesics,}
Bulletin (New Series) of the American Mathematical Society {\bf 35:3} (1998), 193--214. 

\bibitem[S2]{S:S2}
P.~ Schmutz. {\em Systoles on Riemann surfaces,} 
Manuscripta Math., {\bf 85} (1994), 429--447. 

\bibitem[ST]{ST:vol}
G. ~Schumacher and S. ~Trapani.
{\em Estimates of Weil-Petersson volumes via effective 
divisors} Comm. Math. Phys. {\bf 222}, No.1 (2001), 1--7.

\bibitem [SS]{S:S}
S.~ Sabourau. {\em Asymptotic bounds for separating systoles on surfaces,}
Commentarii Mathematici Helvetici, {\bf 83} (2008), no. 1, 35-54.

\bibitem [T]{T}
L. Teo, 
{\em The Weil-Petersson Geometry of the Moduli Space of Riemann Surfaces,}
 Proc. Amer. Math. Soc. {\bf 137} (2009) 541-552. 


\bibitem[Wi]{W} 
E.~Witten. {\em Two-dimensional gravity and intersection theory on moduli spaces,}
 Surveys in Differential Geometry {\bf 1} (1991), 243--269.

 \bibitem[W1]{W:F}
 S.~Wolpert. {\em An elementary formula for the Fenchel-Nielsen twist,}
 Comment. Math. Helv. {\bf 56} (1981), 132--135.

 \bibitem[W2]{W:S}
 S.~Wolpert.
{\em On the symplectic geometry of deformations of a hyperbolic surface,}
 Ann. of Math. {\bf 117:2} (1983), 207--234.

 \bibitem[W3]{W:L}
S.~Wolpert. {\em Behavior of geodesic-length functions on Teichm\"uller space,} 
J. Differential Geom. {\bf 79:2} (2008), 277--334. 

\bibitem[W3]{W:M}
S.~Wolpert. {\em The Weil-Petersson metric geometry,}
 In Handbook of Teichm\"uller theory. Vol. II, volume 13 of IRMA Lect. Math. Theor. Phys., 
47 --64. Eur. Math. Soc., Zurich, 2009.

%\bibitem[Ya]{Yau}
%S. T. Yau.
%{\em Isoperimetric constants and the first eigenvalue of a compact Riemannian manifold,}
 %Ann. Sci. Ecole Norm. Sup. 8 (1975), 487--507.

\bibitem[Z1]{Z:p}
P.~Zograf. The Weil-Petersson volume of the moduli space of punctured spheres, Mapping 
class groups and moduli spaces of Riemann surfaces. Contemp. Math., vol. 150, Amer. Math.
Soc., 1993, 367--372.

\bibitem[Z2]{Z:con}
P.~Zograf. {\em On the large genus asymptotics of Weil-Petersson volumes,} 
 Preprint. 


\end{thebibliography}
\end{document}